\documentclass[12pt,leqno,twoside]{article} 
\usepackage{bezier,amsmath,amssymb,stmaryrd}

\input xy
\xyoption{all}
\input{rk.sty}
\DeclareMathOperator{\ind-crown}{ind-crown}
\DeclareMathOperator{\pro-crown}{pro-crown}
\DeclareMathOperator{\Vm}{V}
\newcommand{\Lm}{\Lambda}
\DeclareMathOperator{\Mor}{Mor}
\newcommand{\ulEl}{\ul{\El\!}\,}
\newcommand{\ulX}{\ul{X\!}\,}

\begin{document}
\title{On lifting stable diagrams in Frobenius categories}
\author{Matthias K\"unzer}
\maketitle

\begin{center}
{\it Dedicated to Claus M.\ Ringel on the occasion of his 60th birthday.} 
\end{center}

\begin{small}
\begin{quote}
\begin{center}{\bf Abstract}\end{center}\vspace*{2mm}
Suppose given a Frobenius category $\El$, i.e.\ an exact category with a big enough subcategory $\Bl$ of bijectives. Let $\ulEl := \El/\Bl$ denote its classical stable category. 
For example, we may take $\El$ to be the category of complexes $\CC(\Al)$ with entries in an additive category $\Al$, in which case $\ulEl$ is the homotopy category of complexes $\KK(\Al)$. 
Suppose given a finite poset $D$ that satisfies the combinatorial condition of being {\it ind-flat}. Then, given a diagram of shape $D$ with values in 
$\ulEl$ (i.e.\ stably commutative), there exists a diagram consisting of pure monomorphisms with values in $\El$ (i.e.\ commutative) that is isomorphic, as
a diagram with values in $\ulEl$, to the given diagram.
\end{quote}
\end{small}

\renewcommand{\thefootnote}{\fnsymbol{footnote}}
\renewcommand{\thefootnote}{\arabic{footnote}}

\begin{footnotesize}
\renewcommand{\baselinestretch}{0.7}
\parskip0.0ex
\tableofcontents
\parskip1.2ex
\renewcommand{\baselinestretch}{1.0}
\end{footnotesize}

\setcounter{section}{-1}

\section{Introduction}

\subsection{The problem}

Let $\El$ be a Frobenius category; that is, an exact category in the sense of {\sc Quillen} \mb{\bfcite{Qu73}{\S 2}} with enough bijective 
objects; cf.\ e.g.\ \bfcite{Ku05}{Sec.\ A.6}. Let $\Bl\tm\El$ denote the full subcategory of bijective objects, and
let $\ulEl = \El/\Bl$ denote the classical stable category of $\El$. Let $\El^{\mono}\tm\El$ denote the subcategory of pure monomorphisms of $\El$. Write $\El\lraa{N}\ulEl$ for the residue class 
functor, and likewise, by abuse of notation, $\El^{\mono}\lraa{N}\ulEl$ for its restriction to $\El^{\mono}$.

Let $D$ be a category. A functor $X$ from $D$ to $\ulEl$ is a diagram of shape $D$ with values in $\ulEl$, sometimes called a {\it stable diagram.} Choosing representatives in $\El$, we may think 
of $X$ as a ``diagram of shape $D$ with values in $\El$ that stably commutes''. We ask under which conditions on $D$ we can find a ``strictly commutative'' diagram $X'$ of shape $D$ with values in 
$\El$ that becomes isomorphic to the ``stably commutative'' diagram $X$, when considering both in the category of diagrams of shape $D$ with values in $\ulEl$.

Put formally, the residue class functor $\El^{\mono}\lraa{N}\;\ulEl$ induces a functor $\El^{\mono}(D)\lrafl{25}{N(D)}\ulEl(D)$ on the diagrams of shape $D$ by pointwise application. We ask for a 
sufficient condition on $D$ for $\El^{\mono}(D)\lrafl{25}{N(D)}\ulEl(D)$ to be dense for all Frobenius categories $\El$; that is, for its induced map on the isoclasses to be surjective.

Such a condition is then a fortiori sufficient for the induced functor $\El(D)\lrafl{25}{N(D)}\ulEl(D)$ to be dense. It turns out to be technically advantageous to consider $\El^{\mono}$ instead of $\El$.

Restricting ourselves to the case of $D$ being a finite poset, we will find a sufficient condition in combinatorial terms on $D$ ensuring that $\El^{\mono}(D)\lrafl{25}{N(D)}\ulEl(D)$ is dense, called 
\mb{\it ind-flatness;} cf.\ Section \ref{SecRes} below.

\subsection{Problems that remain open}

\subsubsection{A precise obstruction against density\,?}

I do not know a necessary and sufficient combinatorial condition on $D$ for $\El^{\mono}(D)\lrafl{25}{N(D)}\ulEl(D)$ to be dense for all Frobenius categories $\El$. 
For instance, it is dense for $D = \De_m\ti\De_n\,$, where $m,\, n\,\ge\, 0$. However, I do not know 
whether it is dense for $D = \De_1\ti\De_1\ti\De_1$.

Considering a category of spaces instead of a Frobenius category $\El$, {\sc Dwyer}, {\sc Kan} and {\sc Smith} have established obstruction classes in certain Hochschild-Mitchell cohomology groups in 
dimension $\ge 3$ against density; cf.\ \bfcite{DKS89}{3.5,\,3.6}.

{\sc Mitchell} gave a combinatorial criterion for the Hochschild-Mitchell cohomology groups to vanish in dimensions $\ge 3$; cf.\ \bfcite{Mi72}{Th.\ 35.7}; cf.\ Section \ref{SecCDKSM}. 
This criterion is fulfilled by $\De_1\ti\De_1$, but not by $\De_1\ti\De_1\ti\De_1$. 

I do not know whether ind-flat finite posets satisfy Mitchell's criterion. I do not know whether there exists an obstruction theory in the spirit of \bfcit{DKS89} for Frobenius categories.
If both should turn out to be true, this would yield the ``true reason'' for density in the case of an ind-flat finite poset. And if, moreover, the obstruction classes should turn out to 
be calculable for $D = \De_1\ti\De_1\ti\De_1$, it would probably also yield an example in which density fails. 

\subsubsection{$1$-Epimorphy\,?}

A functor $\Ul\llaa{F} \Vl$ whose induced functor $\Cl(\Ul)\lrafl{25}{\Cl(F)}\Cl(\Vl)$ given by restriction along $F$ is full and faithful for all categories $\Cl$ is called {\it $1$-epimorphic}; 
cf.\ \bfcite{Ku05}{Sec.\ A.8}. If the finite poset $D$ is a finite quasitree in the sense of Definition \ref{Def4.0.1}, then $\El^{\mono}(D)\lrafl{25}{N(D)}\ulEl(D)$ is \mb{$1$-epimorphic}; see 
Proposition \ref{Th4.1}. We do not know any less drastically restrictive sufficient condition on $D$ for this $1$-epimorphy to hold.

\subsection{Motivation}

The functor $\El^{\mono}(\De_1)\lrafl{25}{N(\De_1)}\ulEl(\De_1)$ being dense can be seen as the technical reason why every morphism in $\ulEl$ can be extended to a distinguished triangle in the sense of 
{\sc Verdier} \bfcit{Ve63}. And the functor $\El^{\mono}(\De_2)\lrafl{25}{N(\De_2)}\ulEl(\De_2)$ being dense can be seen as the main technical reason why the octahedral axiom (TR 4) of loc.\ cit.\ holds. 
We attempt to extend this density property as far as possible.

{\sc Heller} asked the density question in a more general setting; cf.\ \bfcite{He88}{p.\ 4; Prop.\ III.3.9 and remark thereafter}. This question also appeared in the 
discussion of the axioms of a triangulated derivator, due to {\sc Grothendieck} and {\sc Maltsiniotis}; cf.\ \bfcite{Ma01}{p.\ 4}; \mb{cf.\ \bfcit{Ke01}, \bfcit{Gr90}.} 

For applications in topology of the solution of an analogous problem for spaces, see \mb{\bfcite{Co78}{Sec.\ 2}}.

\subsection{Result}
\label{SecRes}

Let $Q$ be a finite poset, considered as a category. For $q\in Q$, let
\[
\barcl
\Lm(q)           & :=  & \{ r\in Q\; :\; r\leq q\} \\
\Lm^0(q)         & :=  & \{ r\in Q\; :\; r < q\} \\
\Vm(q)           & :=  & \{ r\in Q\; :\; r\geq q\} \\
\max(Q)          & :=  & \{ r\in Q \; :\; \Vm(r) = \{ r\} \} \\
\Ob\ind-crown(Q) & := & \Cup_{r,\, s\,\in\,\max(Q)}\Ob \max(\Lm(r)\cap\Lm(s))\; , \\
\ea
\]
yielding a poset $\ind-crown(Q)$ via 
\[
r <_{\ind-crown(Q)} s \Icm:\equ\Icm \mb{$r <_Q s\;\;$ and $\;\; r\not\in\max(Q)\;\;$ and $\;\; s\in\max(Q)\;$}\; .
\] 
We sketch a finite poset $Q$ and its ind-crown.
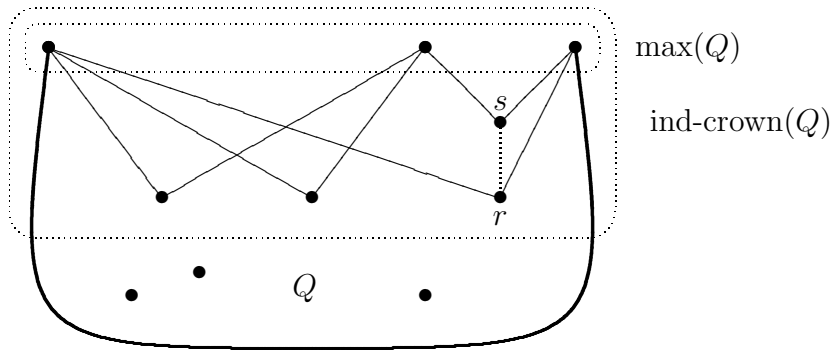
\[
\xy 0;<1cm,0cm>:  
(0,2)*={\bt}; (1.5,0)*={\bt} **@{-}, 
(0,2)*={\bt}; (3.5,0)*={\bt} **@{-} , 
(0,2)*={\bt}; (6,0)*{\bt}*!U(2.9)={r} **@{-} , 
(1.5,0); (5,2)*={\bt} **@{-} , 
(3.5,0); (5,2)*={\bt} **@{-} ,
(6,0); (7,2)*={\bt} **@{-} ,
(5,2); (6,1)*={\bt}*!D(2.9)={s} **@{-} ,
(6,1); (7,2)*={\bt} **@{-} ,
(6,0); (6,1) **@{.} ,
(0,2) *+++{} ; (7,2)*+++{} **\frm<2mm>{.} ,
(8.5,2)*={\max(Q)} ,
(0,2) *+++++{} ; (7,0)*+++++{} **\frm<4mm>{.} ,
(9.2,1)*={\ind-crown(Q)} ,
(0,2);(3.5,-2)**\crv{~*={.}(-0.5,-2)&(-0.5,-2)} ,
(3.5,-2);(7,2)**\crv{~*={.}(7.5,-2)&(7.5,-2)} ,
(3.4,-1.2)*{Q} ,
(5.0,-1.3)*{\bt} ,
(2.0,-1.0)*{\bt} ,
(1.1,-1.3)*{\bt} ,
\endxy     
\]
Whereas it might be the case that $r < s$ in $Q$, we have $r \not< s$ in $\ind-crown(Q)$ since $s\not\in\max(Q)$.

A finite poset $P$ is called {\it ind-flat} if $\ind-crown(\Lm^{\! 0}(p))$ is componentwise $1$-connected for each $p\in P$; cf.\ Definition \ref{Def1.1.2}. For some examples, see 
Definition \ref{Def2.0.1} and Example \ref{Ex2.0.2}.

{\bf Theorem} (Theorem \ref{Th3.1.1}). {\it Suppose given an ind-flat finite poset $D$ and a Frobenius category $\El$. Then $\El^{\mono}(D)\lrafl{25}{N(D)}\ulEl(D)$ is dense.} 

\subsection{Acknowledgement}

I thank the referee for pointing out the work of {\sc Dwyer,} {\sc Kan} and {\sc Smith} \bfcit{DKS89}.

\subsection{Notation and conventions}

\begin{footnotesize}
\begin{itemize}
\item[(i)] For $a,\,b\,\in\,\Z$, we denote by $[a,b] := \{z\in\Z\;:\; a\leq z\leq b\}$ the integral interval. 
\item[(ii)] Given $n\ge 0$, we let $\De_n$ be the linearly ordered set $[0,n]$, with ordering inherited from $\Z$.
\item[(iii)] Given a set $M$, we denote by $\Pfk(M) = \{ N : N\tm M\}$ its power set. If $M$ is finite, then $\# M$ denotes the cardinality of $M$. 
\item[(iv)] All categories are supposed to be small with respect to a sufficiently big universe.
\item[(v)] Composition of morphisms is written on the right, $\lraa{a}\lraa{b} = \lraa{ab}$. 
\item[(vi)] The category of functors and transformations from a category $D$ to a category $\Cl$ is denoted by $\fbo D,\Cl\fbc$, or by $\Cl(D)$. The 
latter is used to emphasise that the objects of $\Cl(D)$ can be viewed as diagrams of shape $D$ with values in $\Cl$; we shall also refer to them as {\it diagrams.}
\item[(vii)] Given a category $\Cl$ and objects $X,\, Y\,\in\,\Ob\Cl$, the set of morphisms from $X$ to $Y$ is denoted by $\liu{\Cl\,}(X,Y)$.
\item[(viii)] Given a category $\Cl$, its opposite category is denoted by $\Cl^\0$.
\item[(ix)] A poset $P = (P,\leq) = (P,\leq_P)$ is a partially ordered set. To consider it as a category, we let $\liu{P}{(p,q)} = \{(p\lra q)\}$ if $p\leq q$, and $\liu{P}{(p,q)} =\leer$ otherwise.
A {\it full subposet} of a poset is a full subcategory. A {\it subposet} is a subcategory. 
\item[(x)] A poset $P$ is {\it discrete} if $p\le q$ implies $p = q$ for $p,\, q\,\in\, P$; that is, if each morphism in $P$ is an identity.
\item[(xi)] Given an exact category $\El$, we denote by $\El^{\mono}$ its subcategory of pure monomorphisms, and by $\El^{\epi}$ its subcategory of pure epimorphisms. 
By $\lramono$, we denote a pure monomorphism; by $\lraepi$, we denote a pure epimorphism. Cf.\ e.g.\ \bfcite{Ku05}{Sec.\ A.2}.
\item[(xii)] A {\it Frobenius category} $\El$ is an exact category in which each $X\in\Ob\El$ allows for $N\lraepi X\lramono N'$ with bijective objects $N$ and $N'$; 
cf.\ e.g.\ \mb{\bfcite{Ku05}{Sec.\ A.2.3}}. Denoting by $\Bl\tm\El$ its full subcategory of bijective objects, we let $\ulEl := \El/\Bl$ denote the classical stable category of $\El$.
Given a morphism $X\lrafl{25}{f} Y$ in $\El$, its residue class in $\ulEl$ is denoted by $\ulX\lrafl{30}{\ul{f\!}\,} \ul{Y\!\!}\,\,$.
\end{itemize}
\end{footnotesize}
\section{Limits and pure monomorphisms}

\subsection{Crowns}

\begin{quote}
\begin{footnotesize}
We extract the relevant part of a poset with respect to taking direct limits of diagrams on it, called its ind-crown, and consider its $1$-connectedness. 
\end{footnotesize}
\end{quote}

\begin{Definition}
\label{Def1.1.1}
\rm
Let $P$ be a finite poset, considered as a category whenever necessary. Given $p\in P$, we define full subposets of $P$
\[
\ba{rcrcl}
\Lm(p) & = & \Lm_P(p)               & := & \{ q\in P\; :\; q\leq p\} \\
\Vm(p) & = & \Vm_{\! P}(p)          & := & \{ q\in P\; :\; q\geq p\} \\
\Lm^{\! 0}(p) & = & \Lm^{\! 0}_P(p) & := & \{ q\in P\; :\; q < p\} \\
\Vm_{\! 0}(p) & = & \Vm_{\! 0,P}(p) & := & \{ q\in P\; :\; q > p\}\; . \\
\ea
\]
Moreover, we define full subposets of $P$
\[
\barcl
\max(P) & := & \{ q\in P \; :\; \Vm(q) = \{ q\} \} \\
\min(P) & := & \{ q\in P \; :\; \Lm(q) = \{ q\} \}\; , \\
\ea
\]
which are discrete. We let
\[
\barcl
\Ob\ind-crown(P) & := & \Cup_{p,\, q\,\in\,\max(P)}\Ob \max(\Lm(p)\cap\Lm(q))  \\
\Ob\pro-crown(P) & := & \Cup_{p,\, q\,\in\,\min(P)}\Ob \min(\Vm(p)\cap\Vm(q))\; .  \\
\ea
\]
The subset $\Ob\ind-crown(P)$ of $\Ob P$ carries a structure of a poset by letting 
\[
p <_{\ind-crown(P)} q \Icm :\equ\Icm \mb{$p <_P q\;$ and $\; p\not\in\max(P)\;$ and $\; q\in\max(P)\;$}
\]
for $p,\, q\,\in\,\Ob\ind-crown(P)$. So $\ind-crown(P)$ is a subposet of $P$, but in general not a full subposet of $P$; cf.\ Example \ref{Ex1.1.5}.

The subset $\Ob\pro-crown(P)$ of $\Ob P$ carries a structure of a poset by letting 
\[
p <_{\pro-crown(P)} q \Icm :\equ\Icm \mb{$p <_P q\;\;$ and $\;\; p\in\min(P)\;\;$ and $\;\; q\not\in\min(P)\; $}
\]
for $p,\, q\,\in\,\Ob\pro-crown(P)$. So $\pro-crown(P)$ is a subposet of $P$, but in general not a full subposet of $P$.

We have $\pro-crown(P) = \ind-crown(P^\0)^\0$.

A poset $C$ is called a {\it crown} if it is finite and if $C = \min(C)\cup \max(C)$. I.e.\ a finite poset $C$ is a crown if there do not exist elements $c,\, c',\, c''\,\in\, C$ with
$c < c' < c''$.

If $P$ is an arbitrary finite poset, then both $\ind-crown(P)$ and $\pro-crown(P)$ are crowns. 
\end{Definition}

\begin{Definition}
\label{Def1.1.2}\rm
Suppose given a crown $C$. Let $\Mor' C$ be the set of non identical morphisms of $C$. 
Let $\Q[\Mor' C]$ be the vector space over $\Q$ with basis $\Mor' C$, and let $\Q[\Ob C]$ be the vector space over $\Q$ with basis $\Ob C$.

The crown $C$ is called {\it componentwise $1$-connected} if the $\Q$-linear map
\[
\barcl
\Q[\Mor' C] & \lraa{\dell_C} & \Q[\Ob C] \\
(c\lra d)   & \lramaps       & d - c     \\
\ea
\]
is injective. Then $C$ is componentwise $1$-connected if and only if $C^\0$ is.
\end{Definition}

\bq
 In other words, a crown $C$ is componentwise $1$-connected if and only if the topological realisation of its nerve is componentwise $1$-connected. In fact, for a finite wedge of circles
 to be $1$-connected, i.e.\ to consist of no circles at all, we may require that $\HH^1$ vanish.
\eq

\begin{Lemma}
\label{Lem1.1.2.1}
If $U\tm C$ is a full subposet of a componentwise $1$-connected crown $C$, then $U$ is itself a componentwise $1$-connected crown.
\end{Lemma}

{\it Proof.}
The poset $U$ is a crown, since there do not exist $c,\,c',\,c''\,\in\, U$ with $c < c' < c''$, for they do not exist in $C$. By restriction, 
injectivity of $\Q[\Mor' C]\lraa{\dell_C}\Q[\Ob C]$ implies injectivity of $\Q[\Mor' U]\lraa{\dell_U}\Q[\Ob U]$.\qed

\begin{samepage}
\begin{Lemma}[recursive characterization]
\label{Lem1.1.3}\Absatz
The crown $C$ is componentwise $1$-connected if and only if {\rm (i)} or {\rm (ii)} or {\rm (iii)} holds.
\begin{itemize}
\item[{\rm (i)}]   There exists $c\in\max(C)$ such that $\#\Lm^{\! 0}(c)\leq 1$, and such that the full subposet $C\ohne\{ c\}$ of $C$ 
                   is componentwise $1$-connected.
\item[{\rm (ii)}]  There exists $c\in\min(C)$ such that $\#\Vm_{\! 0}(c)\leq 1$, and such that the full subposet $C\ohne\{ c\}$ of $C$ 
                   is componentwise $1$-connected.
\item[{\rm (iii)}] $C = \leer$.
\end{itemize}
\end{Lemma}
\end{samepage}

{\it Proof.}
Suppose $C\neq\leer$ to be componentwise $1$-connected. We claim that (i) or (ii) holds.

A {\it chain} in $C$ is a tuple $(c_1,\dots,c_m)$ for some $m\geq 1$ 
such that $c_i < c_{i+1}$ or $c_i > c_{i+1}$ for all $i\in [1,m-1]$, and such that $c_{i+2} \neq c_i$ for all $i\in [1,m-2]$. Suppose given such a chain in $C$.

Assume that there are $j,\, k\,\in\, [1,m]$ such that $j < k$, but $c_j = c_k$. Choose $k - j$ to be minimal with this property. Hence in $(c_j,c_{j+1},\dots,c_{k-1})$, we have
pairwise different entries. The number $k - j$ is even and $\geq 4$. 

If $c_j < c_{j+1}$, then we let 
\[
\ga \; :=\; \sum_{i\in [1,(k-j)/2]} \big((c_{j+2i-2}\lra c_{j+2i-1}) - (c_{j+2i}\lra c_{j+2i-1})\big) \;\in\;\Q[\Mor' C]\; , 
\]
if $c_j > c_{j+1}$, then we let
\[
\ga \; :=\; \sum_{i\in [1,(k-j)/2]} \big((c_{j+2i-1}\lra c_{j+2i-2}) - (c_{j+2i-1}\lra c_{j+2i})\big) \;\in\;\Q[\Mor' C]\; . 
\]
In both cases we have $\ga\neq 0$ since the coefficient of $(c_j\lra c_{j+1})$ resp.\ of $(c_{j+1}\lra c_j)$ equals $1$. In fact, since $c_{j+1}\neq c_{k-1}$, no cancellation occurs. 
But $\ga\dell_C = 0$, and this contradicts the componentwise $1$-connectedness of $C$. From this contradiction we conclude that each chain in $C$ consists of pairwise different entries.

Since $C$ is finite and nonempty, there exists a chain $(c_1,\dots,c_m)$ of maximal length $m$ in $C$. Let $c := c_m$. If $m = 1$, then $c$ satisfies both (i) and (ii). So we may 
suppose $m\geq 2$. We claim that $c$ satisfies (i) if $c_{m-1} < c_m$, and that $c$ satisfies (ii) if $c_{m-1} > c_m$. 

Suppose $c_{m-1} < c_m$. Assume $\#\Lm^{\! 0}(c) > 1$, and let $c_{m+1} \in \Lm^{\! 0}(c)\ohne\{ c_{m-1}\}$. Then $(c_1,\dots,c_{m-1},c_m,c_{m+1})$ is a chain, contradicting
the maximality of $m$. Thus $\#\Lm^{\! 0}(c) \leq 1$. Moreover, $C\ohne\{c\}$ is itself componentwise $1$-connected by Lemma \ref{Lem1.1.2.1}.

Suppose $c_{m-1} > c_m$. Assume $\#\Vm_{\! 0}(c) > 1$, and let $c_{m+1} \in \Vm_{\! 0}(c)\ohne\{ c_{m-1}\}$. Then $(c_1,\dots,c_{m-1},c_m,c_{m+1})$ is a chain, contradicting
the maximality of $m$. Thus $\#\Vm_{\! 0}(c) \leq 1$. Moreover, $C\ohne\{c\}$ is itself componentwise $1$-connected by Lemma \ref{Lem1.1.2.1}.

Conversely, suppose that (i) or (ii) or (iii) holds. We have to show that $C$ is componentwise $1$-connected. By duality, we may assume that (i) holds. 

If $\Lm^{\! 0}(c)\neq\leer$, we write $\Lm^{\! 0}(c) = \{ d\}$. Then the linear map $\Q[\Mor' C]\lraa{\dell_C}\Q[\Ob C]$ decomposes into
\[
\Q[\Mor' (C\ohne\{ c\})]\ds \Q[\{(d\lra c)\}] \;\;\mramh{\smatzz{\dell_{C\ohne\{ c\}}}{0\ru{-2}}{-\w d}{\w c}}\;\; \Q[\Ob(C\ohne\{ c\})]\ds \Q[\{ c\}]\; ,
\]
where we denote by $\w d$ the map that sends $(d\lra c)$ to $d$, and by $\w c$ the map that sends $(d\lra c)$ to $c$. 

If $\Lm^{\! 0}(c) = \leer$, then the linear map $\Q[\Mor' C]\lraa{\dell_C}\Q[\Ob C]$ decomposes as 
\[
\Q[\Mor' (C\ohne\{c\})]\;\;\mra{\smatez{\dell_{C\ohne\{ c\}}}{\; 0}}\;\; \Q[\Ob(C\ohne\{ c\})]\ds \Q[\{ c\}]\; .
\]
In both cases, injectivity of $\dell_C$ results from injectivity of $\dell_{C\ohne\{c\}}$.\qed

\begin{quote}
 \begin{footnotesize}
 \begin{Example}
 \label{Ex1.1.4}\rm
 Let $P = \Pfk(\{ 1,2,3\}) \ohne\big\{\{ 1,2,3\}\big\}$, ordered by inclusion. We have $\max(P) =  \big\{ \{1,2\}, \{1,3\}, \{2,3\} \big\}$. Moreover, we have 
 $\max\!\big(\Lm(\{1,2\})\cap \Lm(\{1,2\})\big) = \big\{\{1,2\}\big\}$, we have $\max\!\big(\Lm(\{1,2\})\cap \Lm(\{2,3\})\big) = \big\{\{2\}\big\}$, etc. Thus,
 \[
 C \; :=\; \ind-crown(P) \= \big\{ \{1,2\}, \{1,3\}, \{2,3\}, \{1\}, \{2\}, \{3\}\big\}\; .
 \] 
 In this example, $C$ is actually a full subposet of $P$. The map $\Q[\Mor' C]\lraa{\dell_C}\Q[\Ob C]$, $(p\lra q)\lramaps q - p$, is given by the matrix
 \[
 \begin{array}{r|rrrrrr}
                     & \{ 1\} & \{ 2\} & \{ 3\} & \{ 1,2\} & \{ 1,3\} & \{ 2,3\} \\\hline
 \{ 1\}\lra\{ 1,2\}  & -1     &  0     &  0     & +1       &  0       &  0       \\
 \{ 1\}\lra\{ 1,3\}  & -1     &  0     &  0     &  0       & +1       &  0       \\
 \{ 2\}\lra\{ 1,2\}  &  0     & -1     &  0     & +1       &  0       &  0       \\
 \{ 2\}\lra\{ 2,3\}  &  0     & -1     &  0     &  0       &  0       & +1       \\
 \{ 3\}\lra\{ 1,3\}  &  0     &  0     & -1     &  0       & +1       &  0       \\
 \{ 3\}\lra\{ 2,3\}  &  0     &  0     & -1     &  0       &  0       & +1       \\
 \end{array}
 \]
 with kernel $\Q\spi{ (+1\; -\!1\; -\!1\; +\!1\; +\!1\; -\!1)}$. Hence the ind-crown $C$ of $P$ is not componentwise $1$-connected. 
 \end{Example} 
 
 \begin{Example}
 \label{Ex1.1.5}\rm
 Let $P = \big\{\, \leer, \{ 1\}, \{ 2\}, \{ 2,3\}, \{2,4\} \big\}$, ordered by inclusion. Then $\Ob\ind-crown(P) = \Ob(P)$. We have $\leer <_P\{2\}$, however, 
 $\leer\not<_{\ind-crown(P)} \{ 2\}$, since $\{ 2\}\not\in\max(P)$. Thus $\ind-crown(P)$ is a subposet of $P$, but not a full subposet.
 Note that $P$ is not a crown, but that, of course, $\ind-crown(P)$ is a crown.
 \end{Example} 

 \begin{Example}
 \label{Ex1.1.6}\rm
 Let $P = \big\{ \{ 1\}, \{ 2\}, \{ 1,2\}, \{ 2,3\} \big\}$, ordered by inclusion. Then $P$ is a crown. We have
 \[
 \ba{rclcl}
 \ind-crown(P) & = & \big\{\{ 2\}, \{ 1,2\}, \{ 2,3\} \big\} & \tmne & P \\
 \pro-crown(P) & = & \big\{\{ 1\}, \{ 2\}, \{ 1,2\} \big\}   & \tmne & P \; .\\
 \ea
 \]
 \end{Example}
\end{footnotesize}
\end{quote}

\subsection{Limits}
\label{SubsecLim}

\begin{quote}
\begin{footnotesize}
We generalise familiar properties of pushouts in exact categories to direct limits over more general diagrams. 
\end{footnotesize}
\end{quote}

Let $\El$ be an exact category; cf.\ e.g.\ \bfcite{Ku05}{Sec.\ A.2}. Let $P$ be a poset. Given a diagram $X\in\Ob\El(P)$, we write
$X(p) =: X_p$ for $p\in\Ob P$, and $X(p\lra q) =: \xi_{p,q}$ whenever $\, p,\, q\,\in\, \Ob P$ with $p \le q$. We write $\limd_P X = \limd_{p\in P} X_p$.
Similarly, the morphisms in a diagram $X'\in\Ob\El(P)$ are denoted by $\xi'_{p,q}$ etc.

\begin{Lemma}
\label{Lem1.2.0_5}
Let $C$ be a componentwise $1$-connected crown, and let $X\in\Ob\El(C)$ be a diagram consisting of pure monomorphisms $\xi_{c,d}$ for all $\, c,\, d\,\in\, C$ with $c\le d$.
Then $\limd_C X\ru{-2.5}$ exists, and the transition morphism $X_c\lra \limd_C X$ is a pure monomorphism for each $c\in C$.
\end{Lemma}

{\it Proof.}
We may assume that $C\neq\leer$. We proceed by induction on $\# C$ and choose $c\in C$ such that condition (i) or (ii) of Lemma \ref{Lem1.1.3} holds. Denote 
$L := \limd_{C\ohne\{ c\}} X|_{C\ohne\{ c\}}$, with transition morphism $X_e\lramonoa{\eta_e} L$ for $e\in C\ohne\{ c\}$.

Consider the case that condition (i) of loc.\ cit.\ holds for $c$. 

If $\Lm^{\! 0}(c) = \leer$, then $\limd_C X = L \ds X_c$, and the transition morphisms are given by 
\mb{$X_e \lramonoa{\smatez{\eta_e}{0}} L\ds X_c$} for $e\in C\ohne\{ c\}$ and by $X_c\lramonoa{\smatez{0}{1}} L\ds X_c$.

If $\Lm^{\! 0}(c)$ consists of one element, say $\Lm^{\! 0}(c) = \{ d\}$, then we consider the pushout
\begin{center}
\begin{picture}(250,250)
\put(   0, 200){$L$}
\put(  50, 210){\vector(1,0){130}}
\put( 110, 210){\ci}
\put( 102, 225){$\scm\lm$}
\put( 200, 200){$\w L$}

\put(  15,  50){\vector(0,1){130}}
\put(  15, 110){\ci}
\put( -35, 110){$\scm\eta_d$}
\put( 215,  50){\vector(0,1){130}}
\put( 215, 110){\ci}
\put( 235, 110){$\scm\mu$}
\put( 170, 170){\line(-1,0){40}}
\put( 170, 170){\line(0,-1){40}}

\put( -10,   0){$X_d$}
\put(  60,  10){\vector(1,0){120}}
\put( 110,  10){\ci}
\put(  90,  28){$\scm\xi_{d,c}$}
\put( 200,   0){$X_c\; $.}
\end{picture}
\end{center}
We have $\limd_C X = \w L$, and the transition morphisms are given by $X_e \lramonoa{\eta_e\lm} \w L$ for $e\in C\ohne\{ c\}$ and by $X_c\lramonoa{\mu} \w L$.

Consider the case that condition (ii) of loc.\ cit.\ holds for $c$. 
We may assume that $\Vm_{\! 0}(c)$ consists of one element, say $\Vm_{\! 0}(c) = \{ d\}$, for otherwise condition (i) of loc.\ cit.\ holds.
We have $\limd_C X = L$, and the transition morphisms are given by $X_e \lramonoa{\eta_e} L$ for $e\in C\ohne\{ c\}$ and by $X_c\;\lramonoa{\xi_{c,d}\eta_d}\; L$.\qed

\begin{quote}
 \begin{footnotesize}
 \begin{Example}
 \label{Ex1.2.3}\rm
 Let $C = \big\{ \{ 1\}, \{ 2\}, \{1,2,3\}, \{ 1,2,4\}\big\}$, ordered by inclusion; the poset $C$ is not componentwise $1$-connected. 
 Denote \mb{$a := \{ 1\}$}, \mb{$b := \{ 2\}$}, \mb{$u := \{ 1,2,3\}$} and \mb{$v := \{1,2,4\}$}. 
 Let $\El = \Z\modl$ be the category of finitely generated $\Z$-modules, with all short exact sequences being pure short exact. Let $X_a = X_b = X_u = X_v = \Z$, let $\xi_{a,u} = 1$, $\xi_{a,v} = 1$,
 $\xi_{b,u} = 1$ and $\xi_{b,v} = m\geq 2$. Then $\limd_C X = \Z/(m-1)$, with transition morphisms $X_u\xrightarrow{1} \Z/(m-1)$ and $X_v\xrightarrow{1} \Z/(m-1)$. 
 The diagram $X$ consists of pure monomorphisms. But none of the transition morphisms to the limit is a pure monomorphism.
 \end{Example}  
 \end{footnotesize}
\end{quote}

\begin{Proposition}
\label{Lem1.2.1}
Suppose given a finite poset $P$ such that $C := \ind-crown(P)$ is componentwise $1$-connected. Suppose given a diagram $X\in\Ob\El(P)$
with $\xi_{p,q}$ purely monomorphic for all $\, p,\, q\,\in\,\Ob P$. The following assertions {\rm (i, ii)} hold.
\begin{itemize}
\item[{\rm (i)}] The limits $\;\limd_C X|_C$ and $\;\limd_P X$ exist in $\El$, and the canonical morphism 
\[
\limd_C X|_C \;\;\lra\;\;\limd_P X
\]
is an isomorphism.
\item[{\rm (ii)}] The transition morphism $X_p\lra \limd_P X$ is a pure monomorphism for $p\in P$.
\end{itemize}
\end{Proposition}

{\it Proof.}
By Lemma \ref{Lem1.2.0_5}, it suffices to prove that, with transition morphisms defined by composition, $L := \limd_C X|_C$ 
is the direct limit of the whole diagram $X$. Denote by $X_c\lramonoa{\eta_c} L$ the transition morphism for $c\in C$.

So for $p\in P$, as transition morphism from $X_p$ to $L$ we take 
\[
(X_p\lramonoa{\tht_p} L) \; :=\; (X_p\lramonoa{\xi_{p,c}} X_c\lramonoa{\eta_c} L) 
\]
for some $c\in\max(P)\tm C$ such that $p\leq c$. We need to show that this definition does not depend on the choice of $c$. So assume given 
$d\in\max(P)\ohne\{ c\}$ such that $p\leq d$. We have to show that $\xi_{p,c}\eta_c = \xi_{p,d}\eta_d$. Note that $p\in \Lm(c)\cap\Lm(d)$. 
Let $e\in\max(\Lm(c)\cap\Lm(d))\tm C$. Then $e\not\in\max(P)$, hence $e <_C c$ and $e <_C d$. Thus we obtain
\[
\xi_{p,c}\eta_c \= \xi_{p,e}\xi_{e,c}\eta_c \= \xi_{p,e}\eta_e \= \xi_{p,e}\xi_{e,d}\eta_d \= \xi_{p,d}\eta_d\; .
\]
As to the universal property of the direct limit, suppose given a family of morphisms $(X_p\lraa{\zeta_p} Z)_{p\in P}$ such that $\xi_{p,q}\zeta_q = \zeta_p$
whenever $\, p,\, q\,\in\, P$ such that $p\leq q$. We obtain an induced morphism $L\lraa{\zeta} Z\ru{5}$ such that $\eta_c\zeta = \zeta_c$ for 
$c\in C$. Uniqueness of $\zeta$ is already given with respect to $C$, so it will hold a fortiori with respect to $P$. It remains to show the existence with respect to $P$, 
that is, it remains to show that $\tht_p\zeta = \zeta_p$ for $p\in P$. In fact, using an element $c\in\max(P)$ with $p\leq c$, we obtain
\[
\tht_p\zeta \= \xi_{p,c}\eta_c\zeta \= \xi_{p,c}\zeta_c \= \zeta_p\; . \vspace*{-3mm}
\]
\qed

\section{Replacement lemmata}

\subsection{Replacement}

\begin{Definition}
\label{Def2.0.1}\rm
A finite poset $D$ is called {\it ind-flat} if $\ind-crown(\Lm^{\! 0}(d))$ is componentwise $1$-connected for each $d\in D$.
Dually, $D$ is called {\it pro-flat} if $\pro-crown(\Vm_{\! 0}(d))$ is componentwise $1$-connected for each $d\in D$.
Altogether, $D$ is called {\it flat} if $D$ is ind-flat and pro-flat.
\end{Definition}

\bq
 \begin{Example}
 \label{Ex2.0.2}\rm\Absit
 \begin{itemize}
 \item[(i)]   The poset $P$ in Example \ref{Ex1.1.4} is ind-flat. It is not pro-flat, since \
              \[
              \pro-crown(\Vm_{\! 0}(\leer)) \= \big\{ \{1\}, \{2\}, \{3\}, \{1,2\}, \{1,3\}, \{2,3\}\big\}
              \]
              is not componentwise $1$-connected.
 \item[(ii)]  The poset $P$ in Example \ref{Ex1.1.5} is flat.
 \item[(iii)] The poset $P$ in Example \ref{Ex1.1.6} is flat.
 \item[(iv)]  The poset $\De_m\ti\De_n$ is flat for $m,\, n\,\geq\, 0$.
 \item[(v)] The poset $\big\{ \leer, \{1\}, \{2\}, \{3\}, \{1,4\}, \{1,5\}, \{1,2,3\}, \{3,4\}, \{3,5\}, \{1,2,3,4,5\} \big\}$ is flat.
 \item[(vi)]   The poset $\De_1\ti\De_1\ti\De_1 \iso \Pfk(\{ 1,2,3\})$ is neither ind-flat nor pro-flat.
 \item[(vii)]  More generally, the poset $\De_1^m \iso \Pfk([1,m])$ is neither ind-flat nor pro-flat for $m\geq 3$.
 \end{itemize}
 \end{Example}

 \begin{Example}
 \label{Ex2.0.1.5}\rm
 If $D$ is a flat finite poset and $D'\tm D$ a full subposet, then $D'$ is not ind-flat in general.

 For instance, let $D = \big\{ \{ 1\}, \{ 2\}, \{1,2\}, \{1,2,3\}, \{1,2,4\}, \{1,2,3,4\} \big\}$, containing the full subposet 
 $D' =  \big\{ \{ 1\}, \{ 2\}, \{1,2,3\}, \{1,2,4\}, \{1,2,3,4\} \big\}$. Then $D$ is flat. In $D'$, however, 
 $\ind-crown\!\big(\Lm_{D'}^{\! 0}(\{1,2,3,4\})\big) = \big\{ \{ 1\}, \{ 2\}, \{1,2,3\}, \{1,2,4\} \big\}$ is not componentwise $1$-connec\-ted, and so $D'$ is not 
 ind-flat.
 \end{Example}
\eq

Suppose given a Frobenius category $\El$; cf.\ e.g.\ \mb{\bfcite{Ku05}{Sec.\ A.2.3}}. Suppose given a finite poset $D$. 

\begin{Definition}
\label{Def2.0.2.5}\rm
A {\it prefunctor} $X$ from $D$ to $\El$ assigns to each object $a$ of $D$ an object $X_a$ of $\El$,
and to each morphism $a\lra b$ of $D$ a morphism $\xi_{a,b}$ of $\El$ in such a way that whenever $a\le b\le c$ 
in $D$, then $\xi_{a,b}\xi_{b,c} - \xi_{a,c}$ is homotopic to zero, i.e.\ it factors over a bijective object in $\El$. Sometimes, we refer to $X$ as a 
{\it prediagram} on $D$ with values in $\El$.

Given prefunctors $X$ and $X'$ from $D$ to $\El$, a morphism $X'\lraa{f} X$ is a tuple $(X'_a\lraa{f_a} X_a)_{a\in\Ob D}$ such that 
$f_a\xi_{a,b} = \xi'_{a,b}f_b$ whenever $a\leq b$ in $D$. Such a morphism $X'\lraa{f} X$ is called a {\it homotopism} if its image 
$\ulX'\lrafl{30}{\ul{f\!}\,} \ulX\ru{5}$ in $\ulEl(D)$ is an isomorphism. 

Let $\El^\sim(D)$ be the category of prefunctors from $D$ to $\El$. In particular, a homotopism is
a morphism in $\El^\sim(D)$. We have a full subcategory $\El(D) \tm \El^\sim(D)$ consisting of diagrams -- a diagram is in particular a prediagram.

There is a canonical dense functor $\El^\sim(D) \lra \ulEl(D)$, $X\lramaps\ulX$, given by taking residue classes of the morphisms that $X$ 
consists of.
\end{Definition}

\begin{Remark}
\label{Lem2.0.3}
Suppose given $X\in\Ob\El^\sim(D)$, a bijective object $N$ in $\El$ and $a\in D$. Let $X'\in\Ob\El^\sim(D)$ be such that
\[
\barcl
X'_b & = & 
\left\{
\ba{ll}
X_b      & \mb{if $b \neq a$} \\
X_a\ds N & \mb{if $b = a$} \\
\ea
\right. \vs\\
(X'_b\;\lrafl{30}{\xi'_{b,c}}\; X'_c) & = & 
\left\{
\ba{ll}
X_b\;\lrafl{26}{\xi_{b,c}}\; X_c                                    & \mb{if $b < c$ and $a\not\in\{ b,c\}$} \vspace*{1mm}\\
X_a\ds N \;\;\lrafl{40}{\smatze{\xi_{a,c}}{\eta_c}}\ru{6.5}\;\; X_c & \mb{if $a = b < c$} \vspace*{1mm} \\
X_b\mra{\smatez{\xi_{b,a}}{\;\zeta_b}} X_a\ds N                     & \mb{if $b < c = a$}\; , \\
\ea
\right. \\
\ea
\]
for some $N\lraa{\eta_c} X_c$ for $c\in\Vm_{\! 0}(a)$ and some $X_b\lraa{\zeta_b} N$ for $b\in\Lm^{\! 0}(a)$. We call $X'$ a {\rm replacement} of $X$ at $a\in D$.

There is an isomorphism 
\[
\ba{rcll}
\ulX'       & \lrafl{18}{\sim}                & \ulX         & \\
X_b         & \lrafl{23}{1}                   & X_b          & \mb{if $b\neq a$} \\
X_a\ds N    & \lrafl{35}{\smatze{1}{0}}\ru{6} & X_a          &  \\
\ea
\]
in $\ulEl(D)$. If $a\in\max(D)$, this isomorphism lifts to a homotopism $X'\lra X$ in $\El^\sim(D)$.
\end{Remark}

\subsection{A purely monomorphic replacement}

\begin{Lemma}
\label{Lem2.1.1}
Suppose given a finite poset $D$ and an element $c\in\max(D)$. Suppose $\ind-crown(\Lm^{\! 0}(c))$ to be componentwise $1$-connected.

Suppose given a diagram $X\in\Ob\El(D)$ such that $X|_{D\ohne\{c\}}\in\Ob\El^{\mono}(D\!\ohne\!\{c\})$, i.e.\ such that its restriction to 
$D\ohne\{ c\}$ consists of pure monomorphisms. Then there exist $X'\in\Ob\El^{\mono}(D)$ and a homotopism $X'\lraa{f} X$.
\end{Lemma}

{\it Proof.}
Let $L := \limd_{\Lm^{\! 0}(c)} X|_{\Lm^{\! 0}(c)}$, which exists in $\El$ by Proposition \ref{Lem1.2.1}.(i). Let $X_b\lramonoa{\eta_b} L$ denote the transition morphism for $b\in \Lm^{\! 0}(c)$, 
which is purely monomorphic by Proposition \ref{Lem1.2.1}.(ii). Let $L\lraa{\zeta} X_c$ be the unique morphism such that $\eta_b\zeta = \xi_{b,c}$ for all $b\in\Lm^{\! 0}(c)$. Choose a pure monomorphism 
$L\lramonoa{\io} N$ with $N$ bijective. For a replacement at $c$ in the sense of Remark \ref{Lem2.0.3}, we let $X'_c := X_c\ds N$ and 
\[
(X'_b\;\lrafl{30}{\xi'_{b,c}}\; X'_c) \; :=\; (X_b\mramonofl{27}{\smatez{\xi_{b,c}}{\;\eta_b\io}} X_c\ds N)
\]
for $b\in\Lm^{\! 0}(c)$. This yields a diagram $X'\in\Ob\El^{\mono}(D)$. Since $c\in\max(D)$, Remark \ref{Lem2.0.3} gives a homotopism $X'\lra X$.
\qed

\begin{Lemma}
\label{Lem2.1.2}
Given a ind-flat finite poset $D$ and a diagram $X\in\Ob\El(D)$. Then there exist $X'\in\Ob\El^{\mono}(D)$ and a homotopism $X'\lraa{f} X$.
\end{Lemma}

{\it Proof.}
We proceed by induction on $\# D$ and may assume $\# D\geq 1$.
Let $c\in\max(D)$. Since $D\ohne\{c\}$ is ind-flat, too, we may assume the assertion to hold for the diagram $X|_{D\ohne\{c\}}$
on $D\ohne\{c\}$; i.e.\ we may assume a homotopism $Y\lraa{g} X|_{D\ohne\{c\}}$ in $\El(D\!\ohne\!\{c\})$ to exist for some $Y\in\Ob\El^{\mono}(D\!\ohne\!\{c\})$. 
Define $X''\in\Ob\El(D)$ by 
\[
\ba{lcll}
X''|_{D\ohne\{c\}}                        & = & Y                                             & \\
X''_c                                     & = & X_c                                           & \\
(X''_b \;\lrafl{30}{\xi'_{b,c}}\; X''_c)  & = & (Y_b\lraa{g_b} X_b \lrafl{30}{\xi_{b,c}} X_c) & \mb{for $b\in D\ohne\{ c\}$}\; .\\
\ea
\]
In $\El(D)$, we have a homotopism $X''\lraa{f} X$ given by
\[
\ba{lcll}
(X''_b \lrafl{30}{f_b} X_b) & = & (Y_b\lrafl{30}{g_b} X_b)         & \mb{for $b\in D\ohne\{ c\}$}\\
(X''_c \lrafl{30}{f_c} X_c) & = & (X_c\lrafl{30}{1_{X_c}} X_c) & \; .\\
\ea
\]
Finally, by Lemma \ref{Lem2.1.1}, we can replace $X''$ by an object $X'$ in $\El^{\mono}(D)$.\qed

\subsection{A replacement that adds a commutativity}

\begin{Lemma}
\label{Lem2.2.1}
Suppose given a finite poset $D$, an element \mb{$c\in\max(D)$}, an element
\mb{$d\in\max(\Lm^{\! 0}(c))$}, and an element \mb{$e\in\Lm^{\! 0}(d)$}. So $e < d < c$, and there is no element in between $d$ and $c$. Suppose 
$\ind-crown(\Lm^{\! 0}(d))$ to be componentwise $1$-connected. 

Suppose given $X\in\Ob\El^\sim(D)$ such that {\rm (I,\;II)} hold.
\begin{itemize}
\item[{\rm (I)}] We have $X|_{D\ohne\{c\}}\in\Ob\El(D\!\ohne\!\{c\})$. 
\item[{\rm (II)}] We have $X|_{\Lm^0(c)}\in\Ob\El^{\mono}(\Lm^0(c))$. 
\end{itemize}
Then there exist $X'\in\Ob\El^\sim(D)$ and an isomorphism $\ulX'\lraiso \ulX$ in $\ulEl(D)$ such that \mb{\rm (i,\;ii,\;iii,\;iv)} hold.
\begin{itemize}
\item[{\rm (i)}] We have $X'|_{D\ohne\{ c\}}\in\Ob\El(D\!\ohne\!\{c\})$.
\item[{\rm (ii)}]  We have $X'|_{\Lm^0(c)}\in\Ob\El^{\mono}(\Lm^0(c))$.
\item[{\rm (iii)}]   We have $\xi'_{e,c} = \xi'_{e,d}\,\xi'_{d,c}\,$.
\item[{\rm (iv)}]  We have $X'|_{D\ohne\{ d\}} \iso X|_{D\ohne\{ d\}}$ in $\El^\sim(D\!\ohne\!\{d\})$.
\end{itemize}
\end{Lemma}

{\it Proof.}
Denote $L := \limd_{\Lm^{\! 0}(d)} X|_{\Lm^{\! 0}(d)}$, and let $X_b\lramonoa{\eta_b} L$ be the transition morphism for $b\in\Lm^{\! 0}(d)$; cf.\ Proposition \ref{Lem1.2.1}.
Let $L\lraa{\zeta} X_d$ be the unique morphism such that $\eta_b\zeta = \xi_{b,d}$ for all $b\in\Lm^{\! 0}(d)$. Choose a pure monomorphism 
$L\lramonoa{\io} N$ with $N$ bijective. Bijectivity of $N$ together with pure monomorphy of $\eta_e\io$ allows to factor the nullhomotopic difference $\xi_{e,c} - \xi_{e,d}\,\xi_{d,c}$ as
\[
\xi_{e,c} - \xi_{e,d}\,\xi_{d,c} \= \eta_e\,\io\,\tht
\]
for some $N\lraa{\tht} X_c\,$.

For a replacement at $d$ in the sense of Remark \ref{Lem2.0.3}, we let $X'_d := X_d\ds N$ and
\[
\ba{rcll}
(X'_b\;\lrafl{30}{\xi'_{b,d}}\; X'_d) & := & (X_b\mramonofl{30}{\smatez{\xi_{b,d}}{\;\eta_b\io}} X_d\ds N) & \mb{for $b\in\Lm^{\! 0}(d)$} \vspace*{3mm}\\
(X'_d\;\lrafl{30}{\xi'_{d,c}}\; X'_c) & := & (X_d\ds N\;\lrafl{45}{\smatze{\xi_{d,c}\ru{-1.5}}{\tht}}\; X_c) & \vspace*{3mm}\\
(X'_d\;\lrafl{30}{\xi'_{d,a}}\; X'_a) & := & (X_d\ds N\;\lrafl{45}{\smatze{\xi_{d,a}\ru{-1.5}}{0}}\; X_a) & \mb{for $a\in\Vm_{\! 0}(d)\ohne\{c\}$}\; . \\
\ea
\]
This yields the required diagram $X'$.\qed

\section{Density}

\begin{Theorem}
\label{Th3.1.1}
Suppose given an ind-flat finite poset $D$. Then the residue class functor
\[
\barcl
\El^{\mono}(D) & \lra     & \ulEl(D) \\
X              & \lramaps & \ulX \\ 
\ea
\]
is dense. 
\end{Theorem}

{\it Proof.}
We proceed by induction on $\# D$. We may assume $\# D\geq 1$. Let $c\in\max(D)$. Suppose given $X\in\Ob\El^\sim(D)$.
Since \mb{$D\ohne\{c\}$} is ind-flat, by induction, there exists a diagram $Y\in\Ob\El^{\mono}(D\ohne\{ c\})$ such that $\ul{Y\!\!}\,\,\lraisoa{g}\ulX$. Extending $Y$ to a diagram 
$\h Y\in\El^\sim(D)$ by appending $X_c$ at $c$, and morphisms $Y_d\lrafl{30}{\h g_d\xi_{d,c}}X_c\ru{5.8}$ for $d < c$, where $\h g_d$ is a representative of $g_d$, 
we obtain $\ul{\h Y\!\!}\,\,\iso\ulX$ via an isomorphism that restricts to $g$ on $D\ohne\{c\}$ and to the identity on $\{c\}$. Moreover, $\h Y|_{D\ohne\{c\}}\in\Ob\El^{\mono}(D)\ru{4.5}$. 
So we may assume that \mb{$X|_{D\ohne\{c\}}\in\Ob\El^{\mono}(D)$}.

A full subposet $U\tm\ind-crown(\Lm^0(c))$ is called {\it commutant} (with respect to $X$) whenever there exist 
\mb{$X'\in\Ob\El^\sim(D)$} and an isomorphism $\ulX'\lraiso\ulX$ such that (1), (2) and (3) hold.
\begin{itemize}
\item[(1)] We have $X'|_{D\ohne\{c\}}\in\Ob\El(D\!\ohne\!\{c\})$. 
\item[(2)] We have $X'|_{\Lm^0(c)}\in\Ob\El^{\mono}(\Lm^0(c))$. 
\item[(3)] We have $\xi'_{s,t}\,\xi'_{t,c} = \xi'_{s,c}$ for all $s,\, t\,\in\, U$ with $s < t$.
\end{itemize}

By assumption, $\ind-crown(\Lm^0(c))$ is componentwise $1$-connected, so by Lemma \ref{Lem1.1.2.1}, any full subposet  
$U\tm\ind-crown(\Lm^0(c))$ is a componentwise $1$-connected crown, too.

We {\it claim} that each full subposet $U\tm\ind-crown(\Lm^0(c))$ is commutant.

We perform an induction on $\# U$. We may assume $\# U\geq 1$. By Lemma \ref{Lem1.1.3}, we can distinguish the following two cases.

Case (i). There exists $u\in\max(U)$ such that $\#\Lm^0_U(u) \leq 1$. 
If $\Lm^0_U(u) = \leer$, then we conclude from $U\ohne\{ u\}$ being commutant that $U$ is 
commutant. So suppose that, say, $\Lm^0_U(u) = \{ v\}$. By induction, we may assume that $\xi_{s,t}\,\xi_{t,c} = \xi_{s,c}$ for all $s,\,t\,\in\, U\ohne\{ u\}$ 
with $s < t$. We use Lemma \ref{Lem2.2.1} in the following way. In the notation used there, we let $c = c$, $d = u$ and $e = v$, and get an $X'\in\Ob\El^\sim(D)$ and an isomorphism 
$\ulX'\lraiso \ulX$ such that $\xi'_{s,t}\,\xi'_{t,c} = \xi'_{s,c}\ru{-2}$ for all $s,\,t\,\in\, U\ohne\{ u\}$ with $s < t$ by loc.\ cit.\ (iv), and 
such that $\xi'_{v,u}\,\xi'_{u,c} = \xi'_{v,c}\ru{-2}$ by loc.\ cit.\ (iii). Finally, $X'|_{D\ohne\{ c\}}\in\Ob\El(D\ohne\{ c\})$ by loc.\ cit.\ (i) and 
$X'|_{\Lm^0(c)}\in\Ob\El^{\mono}(\Lm^0(c))$ by loc.\ cit.\ (ii). Thus $U$ is commutant.

Case (ii). There exists $u\in\min(U)$ such that $\#\Vm_{0,U}(u) \leq 1$. If $\Vm_{0,U}(u) = \leer$, then we 
conclude from $U\ohne\{ u\}$ being commutant that $U$ is commutant. So suppose that, say, $\Vm_{0,U}(u) = \{ v\}$. By induction, we may assume that
$\xi_{s,t}\,\xi_{t,c} = \xi_{s,c}$ for all $s,\,t\,\in\, U\ohne\{ u\}$ with $s < t$. We define $X'\in\Ob\El^\sim(D)$ by letting $\xi'_{s,t} := \xi_{s,t}$ if 
$s,\, t\,\in\, D$ with $s < t$ and $(s,t)\neq (u,c)$, and letting $\xi'_{u,c} := \xi_{u,v}\,\xi_{v,c} = \xi'_{u,v}\,\xi'_{v,c}$. 
Then $\ulX' = \ulX$ and \mb{$\xi'_{s,t}\,\xi'_{t,c} = \xi'_{s,c}\ru{-2}$} for all $s,\,t\,\in\, U$ with \mb{$s < t$.}
Moreover, $X'|_{D\ohne\{ c\}} = X|_{D\ohne\{ c\}}\in\Ob\El(D\ohne\{c\})$ and \mb{$X'|_{\Lm^0(c)} = X|_{\Lm^0(c)}\in\Ob\El^{\mono}(\Lm^0(c))$}. Thus $U$ is commutant.

This proves the {\it claim.} In particular, $\ind-crown(\Lm^0(c))$ is commutant, and we dispose of an according diagram $X'\in\Ob\El^\sim(D)$ satisfying (1), 
(2) and (3).

Now define $X''\in\Ob\El^\sim(D)$ by letting $\xi''_{b,d} := \xi'_{b,d}\ru{-2}$ for $b < d\neq c$ and $\xi''_{b,c} := \xi'_{b,t}\,\xi'_{t,c}$ for 
\mb{$b\in\Lm^0(c)$,} for some $t\in\max(\Lm^0(c))$ with $b\leq t$. Since $\xi'_{s,t}\,\xi'_{t,c} = \xi'_{s,c}\ru{-2}$ for all 
\mb{$s,\,t\,\in\,\ind-crown(\Lm^0(c))$} with $s < t$, this definition of $\xi'_{b,c}$ does not depend on the choice of $t$, and we have in fact $X''\in\Ob\El(D)$ with $\ulX'' = \ulX'$.

By Lemma \ref{Lem2.1.2}, there exist $X'''\in\Ob\El^{\mono}(D)$ and a homotopism $X'''\lra X''$.\qed

\begin{Scholium}
\label{Sch3.1.3}
Given a flat finite poset $D$, the residue class functors $\El^{\mono}(D) \lra \ulEl(D)$ and $\El^{\epi}(D)\lra\ulEl(D)$ are dense.
\end{Scholium}

\bq
\begin{Example}
\label{Ex3.1.3.5}\rm
{\it Given $X\in\Ob\El^\sim(D)$, in general there does not exist $X'\in\Ob\El(D)$ and a homotopism $X'\lra X$.}

Given a finite poset $D$ such that $D\ti\De_1$ is ind-flat, this failure prevents us from using density of 
$\El^{\mono}(D\ti\De_1)\lra\ulEl(D\ti\De_1)$ together with \bfcite{Ku05}{Lem.\ A.35} to conclude that $\El^{\mono}(D)\lra \ulEl(D)$ is $1$-epimorphic.

{\it Proof.} Let $D = \De_2$. Let $\El$ be a Frobenius category in which not every object is bijective. Let $X\in\El^\sim(D)$ be defined to have a non-bijective object 
$X_0$, an arbitrary object $X_1$ and a bijective object $X_2$ such that there exist $X_0\lramonoa{i} X_2\ru{5.5}$; and by morphisms 
$\xi_{0,1} = 0$, $\xi_{1,2} = 0$ and $\xi_{0,2} = i$. 

Assume there is a homotopism $X'\lra X$ for some $X'\in\Ob\El(D)$, consisting of morphisms $X'_i\lraa{u_i} X_i$ for $i\in [0,2]$. 
Then $u_1 \xi_{1,2} = \xi'_{1,2} u_2$ shows that $\xi'_{1,2} u_2 = 0$. Hence 
\[
u_0 i \= u_0 \xi_{0,2} \= \xi'_{0,2} u_2 \= \xi'_{0,1}\xi'_{1,2} u_2 \= 0\; .
\]
Since $i$ is monomorphic, this implies $u_0 = 0$. Since $\ul{u_0}$ is an isomorphism, we conclude that $\ul{X_0}\iso 0$, i.e.\ that $X_0$ is bijective,
contradicting our assumption. Thus there does not exist a homotopism $X'\lra X$ with $X'\in\Ob\El(D)$. \qed
\end{Example}

\begin{Question}
\label{Qu3.1.4}\rm
Is there a poset $D$ and a Frobenius category $\El$ such that the residue class functor 
$\El^{\mono}(D)\lra \ulEl(D)$ is {\bf not} dense\,? What about, say, $D = \De_1\ti\De_1\ti\De_1$\,? 
Is there a counterexample if we relax the condition on $D$ and allow $D$ to be an arbitrary finite category\,?
\end{Question}

To illustrate the kind of problem addressed in Question \ref{Qu3.1.4}, we briefly report a failed attempt to find a counterexample.

\begin{Example}
\label{Qu3.1.5}\rm
We let the finite category $D$ defined by $\Ob D = \{ c\}$ and by $\liu{D}{(c,c)} = \{1_c,\,\al\}$, where $\al\neq 1_c\,$, but $\al^2 = 1_c\,$. 
Let $X := (C\lraa{a} C)$ be an endomorphism of $\El$ that is an object of $\El^\sim(D)$, i.e.\ assume $a^2 - 1$ to vanish in $\ulEl$. Let $C\lramonoa{u} N$ be a pure monomorphism into a bijective object.
Consider a factorization $a^2 - 1 = uv$ and a prolongation $N\lraa{\w a} N\ru{4}$ of $a$ along $u$, i.e.\ $u\w a = a u$. Note that $u (\w a v - v a) = 0$ and $u(\w a^2 - 1 - vu) = 0$. 
 
Assume that $u$, $v$ and $\w a$ can be chosen such that the following hold.
\begin{itemize}
\item[(1)] We have $\w a v - v a = 0$.
\item[(2)] We have $\w a^2 - 1 - vu = 0$. 
\end{itemize}
E.g.\ we might take $\El = \Z/27\modl$, $C = \Z/9$, $a = 2$, $N = \Z/27$, $u = 3$, $v = 1$ and $\w a = 2$.

Let $X'\in\Ob\El(D)$ be defined by 
$
C\ds N\;\,\lrafl{35}{\rsmatzz{a}{u}{-v}{-\w a}\;}\; C\ds N \ru{5.5}
$.
Then $\ulX \iso \ulX'$ in $\ulEl(D)$ via $C\;\lraa{\smatez{1}{0}}\;C\ds N$. So in order to find a counterexample in this manner, it is necessary to use an endomorphism $a$ for which, for all 
choices of $v$ and $\w a$, condition (1) or (2) fails.
\end{Example}
\eq
\section{$1$-Epimorphy}

\begin{Definition}
\label{Def4.0.1}\rm
A finite poset $D$ is called a {\it quasitree} if for all $a,\,b\,\in\, D$, the full subposet $\Vm_0(a)\cap\Lm^0(b)$ of $D$ is linearly ordered.
\end{Definition}

\bq
\begin{Example}
\label{Ex4.0.0.5}\rm
Suppose given a finite poset $D$.
\begin{itemize}
\item[(i)] If $D$ is a crown, then it is a quasitree, since then $\Vm_0(a)\cap\Lm^0(b) = \leer$ for all $a,\,b\,\in\, D$.
\item[(ii)] If for $a,\, b\,\in\, D$ such that $a\not\le b$ and $a\not\ge b$, we have $\Vm(a)\cap\Vm(b) = \leer$, then the poset $D$ is called an {\it ascending tree.} An ascending tree is a quasitree.
\item[(iii)] The poset $D$ is a quasitree if and only if its full subposet $\Vm(a)$ is an ascending tree for all $a\in D$.
\end{itemize}
\end{Example}
\eq

\begin{Lemma}
\label{Lem4.0.2}
Suppose given a finite poset $D$. The following are equivalent.
\begin{itemize}
\item[{\rm (i)}] The poset $D$ is a finite quasitree.
\item[{\rm (ii)}] The subposet $\ind-crown(\Lm^0(a))$ of $D$ is discrete for all $a\in D$.
\item[{\rm (iii)}] The subposet $\pro-crown(\Vm^0(a))$ of $D$ is discrete for all $a\in D$.
\end{itemize}
In particular, if $D$ is a finite quasitree, then $D$ is flat.
\end{Lemma}

{\it Proof.}
First of all, we remark that $\ind-crown(\Lm^0(a))$ is discrete if and only if $\ind-crown(\Lm^0(a)) = \max(\Lm^0(a))$, i.e.\ if and only if
\[
\Lm^0(b)\cap\Lm^0(b') \= \leer
\]
for all $b,\,b'\,\in\,\max(\Lm^0(a))$ with $b\neq b'$.

Ad (i) $\imp$ (ii). Suppose given $b,\,b'\,\in\,\max(\Lm^0(a))$ with $b\neq b'$. Assume there exists $c\in\Lm^0(b)\cap\Lm^0(b')$.
Then $b,\, b'\,\in\,\Vm_0(c)\cap \Lm^0(a)$, but $b\not\le b'$ and $b\not\ge b'$ because of their maximality in $\Lm^0(a)$. But $\Vm_0(c)\cap \Lm^0(a)$
is linearly ordered. This contradicion shows that $\Lm^0(b)\cap\Lm^0(b') = \leer$.

Ad (ii) $\imp$ (i). Given $a,\, c\,\in\, D$, we have to show that $\Vm_0(c)\cap \Lm^0(a)$ is linearly ordered. Assume there exist
$b$ and $b'$ in $\Vm_0(c)\cap\Lm^0(a)$ such that $b\not\le b'$ and $b\not\ge b'$. Choose \mb{$d\in\min(\Lm(a)\cap \Vm_0(b)\cap\Vm_0(b'))$}. Choose \mb{$e\in\max(\Vm(b)\cap\Lm^0(d))$}. 
Choose  \mb{$e'\in\max(\Vm(b')\cap\Lm^0(d))$}. Then $e$ and $e'$ are different elements of $\max(\Lm^0(d))$, because $e = e'$ would imply 
$e \not\in\{b,b'\}$, and we could replace $d$ by $e$, contradicting the minimality of $d$. We have $e,\, e'\in \max(\Lm^0(d))$, whereas
\[
c\;\in\; \Lm^0(e)\cap\Lm^0(e')\;\neq\;\leer\; ,
\]
which is impossible by (ii). This contradiction shows that $\Vm_0(c)\cap\Lm^0(a)$ is in fact linearly ordered.\qed

A functor $\Ul\llaa{F}\Vl$ is called {\it $1$-epimorphic} if the induced functor $\Cl(\Ul)\lrafl{25}{\Cl(F)}\Cl(\Vl)$, given by restriction along $F$, is full and faithful for any category $\Cl$;
cf.\ \bfcite{Ku05}{Sec.\ A.8}.

\begin{Proposition}
\label{Th4.1}
Suppose given a finite quasitree $D$. Then the residue class functor
\[
\barcl
\El^{\mono}(D) & \lra     & \ulEl(D) \\
X              & \lramaps & \ulX \\ 
\ea
\]
is $1$-epimorphic. 
\end{Proposition}

{\it Proof.}
By Lemma \ref{Lem4.0.2} and Theorem \ref{Th3.1.1}, this functor is dense. So by \bfcite{Ku05}{Lem.\ A.35}, it suffices to show that for $X,\, Y\,\in\,\Ob\El^{\mono}(D)$ and a 
morphism $\ulX \lraa{f} \ul{Y\!\!}\,\,$, there exist a homotopism $X'\lraa{g'} X$ and a morphism $X'\lraa{g} Y$ in $\El(D)$ such that
\[
(\ulX'\lrafl{30}{\ul{g}'} \ulX \lrafl{30}{f} \ul{Y\!\!}\,\,) \= (\ulX'\lrafl{30}{\ul{g}} \ul{Y\!\!}\,\, ) \; .
\]
The morphisms that $X$ consists of are denoted by $\xi_{a,b}\,$, the morphisms that $Y$ consists of by $\eta_{a,b}$\,, etc., where $a,\, b\,\in\, D$ with $a < b$.

We proceed by induction on $\# D$. We may assume $\# D\ge 1$. Let $c\in\max(D)$. By induction, the assertion holds for $D\ohne\{ c\}$.
Letting $X''_c := X_c$, by composition, we obtain a diagram $X''\in\Ob\El(D)$, a homotopism $X''\lraa{h'} X\ru{5}$ and a morphism $X''|_{D\ohne\{c\}}\lraa{h} Y|_{D\ohne\{c\}}$ such that the 
following hold.
\begin{itemize}
\item[(i)] The diagram $X''|_{D\ohne\{c\}}$ is in $\Ob\El^{\mono}(D\ohne\{ c\})$.
\item[(ii)] We have $\;\ul{h\!}\,'|_{D\ohne\{c\}} f|_{D\ohne\{c\}} \= \ul{h\!}\,\;$.
\item[(iii)] We have $h'_c = 1_{X_c}$.
\end{itemize}
We choose a representative $X''_c\lraa{{\h f}_c} Y_c$ in $\El$ of $f_c$. We choose a pure monomorphism
\[
\Ds_{a\,\in\,\max(\Lm^0(c))} X_a\;\;\mramono{(i_a)_a}\;\; N
\]
into a bijective object $N$. In particular, each $i_a$ is purely monomorphic. We have a factorisation
\begin{center}
\begin{picture}(450,250)
\put(-120, 200){$\Ds_{a\,\in\,\max(\Lm^0(c))} X_a$}
\put( 200, 210){\vector(1,0){180}}
\put( 280, 210){\ci}
\put( 260, 225){$\scm (i_a)_a$}
\put( 400, 200){$N$}
\put(  15, 170){\vector(0,-1){120}}
\put(-220, 105){$\scm (h_a\eta_{a,c} - \xi''_{a,c}{\h f}_c )_a$}
\put( 380, 190){\vector(-2,-1){325}}
\put( 200, 120){$\scm s$}
\put(   0,   0){$Y_c$}
\end{picture}
\end{center}

Define a replacement $X'$ of $X''$ at $c$ in the sense of Remark \ref{Lem2.0.3} by $X'_c := X''_c\ds N$ and by
\[
(X'_b \lrafl{30}{\xi'_{b,c}} X'_c) \; := (X''_b \;\;\mramonofl{27}{\smatez{\xi''_{b,a}\xi''_{a,c}\;}{\xi''_{b,a}i_a}}\;\; X_c\ds N)\; 
\]
for $b\in\Lm^0(c)$, where $\{ a\} = \max(\Vm(b)\cap\Lm^0(c))$, which is welldefined since $D$ is a quasitree. Then $X'\in\Ob\El^{\mono}(D)$. Let 
$X'\lraa{h''} X''\ru{5}$ be the homotopism of Remark \ref{Lem2.0.3}, and let $(X'\lraa{g'} X) := (X'\lraa{h''} X''\lraa{h'} X)\ru{5}$.
Let $X'\lraa{g} Y$ be defined by
\[
\left\{
\ba{rcll}
(X'_b\lraa{g_b} Y_b) & := & (X''_b\lraa{h_b} Y_b)                          & \mb{at $b\neq c$} \vspace*{4mm}\\
(X'_c\lraa{g_c} Y_c) & := & (X_c\ds N\lrafl{40}{\smatze{{\h f}_c}{s}} Y_c) & \mb{at $c$}\; .\\
\ea
\right.
\]
We claim that $\ul{g\!}\,' f = \ul{g\!}\,$. If $b\neq c$, we obtain
\[
(\ul{g\!}\,' f)_b \= \ul{h\!}\,'_b f_b \= \ul{h\!}\,_b \= \ul{g\!}\,_b \; .
\]
At $c$, we obtain
\[
(\ul{g\!}\,' f)_c \= \smatze{1}{0} f_c \= \smatze{f_c\ru{-2}}{0} \= \ul{\smatze{{\h f}_c}{s}} \= \ul{g\!}\,_c \; .\vspace*{-6mm}
\]
\qed

\begin{Scholium}
\label{Sch4.3}
Given a finite quasitree $D$, the residue class functors $\El^{\mono}(D) \lra \ulEl(D)$ and $\El^{\epi}(D)\lra\ulEl(D)$ are dense and $1$-epimorphic.

\rm
Using Lemma \ref{Lem4.0.2}, this summarises Scholium \ref{Sch3.1.3}, Proposition \ref{Th4.1} and its dual assertion in the given situation.
\end{Scholium}

\bq
 \begin{Example}
 \label{Ex4.5}\rm
 {\it Given a finite quasitree $D$, the residue class functor $\El^{\mono}(D)\lra\ulEl(D)$ is not full in general.}

 A full and dense functor is $1$-epimorphic; cf.\ \bfcite{Ku05}{Cor.\ A.37}. This example, together with Scholium \ref{Sch4.3}, shows that this implication is strict. 

 {\it Proof.} Let $D = \big\{ \{1\},\, \{2\},\,\{1,2\}\big\}$, let $p\ge 3$ be a prime, and let $\El = \Z/p^3\modl$, with all short exact sequences being purely short exact. 
 An object is bijective if and only if it is a finite direct sum of copies of $\Z/p^3$. Consider the following morphism in $\ulEl(D)$.
 \begin{center}
 \begin{picture}(350,350)
 \put( 100, 300){$\Z/p^2$}
 \put( 150, 280){\vector(1,-1){50}}
 \put( 172, 258){\ci}
 \put( 180, 275){$\smatez{p-1\;\;\;}{1}$}
 \put( 200, 200){$\Z/p^2\dk\Z/p^2$}
 \put( -20, 200){$\Z/p^2$}
 \put(  65, 210){\vector(1,0){125}}
 \put( 115, 210){\ci}
 \put(  55, 170){$\smatez{p+1\;}{-1}$} 
 
 \put( 130, 280){\line(0,-1){60}}
 \put( 130, 200){\vector(0,-1){60}}
 \put( 260, 180){\vector(0,-1){140}}
 \put( 270, 110){$\smatze{1\ru{-2}}{1}$}
 \put(  20, 180){\vector(0,-1){140}}

 \put( 122,  95){$0$}
 \put( 150,  80){\vector(1,-1){50}}
 \put( 172,  58){\ci}
 \put( 222,   0){$\Z/p^2$}
 \put(  12,   0){$0$}
 \put(  50,  10){\vector(1,0){140}}
 \put( 110,  10){\ci}
 \end{picture}
 \end{center}

 The question whether it lifts to a morphism in $\El^{\mono}(D)$ is equivalent to the question whether there exist $h,\, k\,\in\,\Z/p$ such that 
 \[
 \barcl
 \smatez{p-1\;\;\;\:}{1}\left(\smatze{1\ru{-2}}{1} + p\smatze{h\ru{-2}}{k}\right) & \con_{p^2} & 0 \vspace*{2mm}\\
 \smatez{p+1\;}{ -1}\left(\smatze{1\ru{-2}}{1} + p\smatze{h\ru{-2}}{k}\right)     & \con_{p^2} & 0 \; .\\
 \ea
 \]
 Adding the two resulting equations, we get $2p\con_{p^2} 0$, so that we cannot find the required $h$ and $k$.\qed
\end{Example}

\begin{Question}
\label{Qu4.4}\rm
Given an ind-flat finite poset $D$ and a Frobenius category $\El$, is the residue class functor $\El^{\mono}(D)\lra \ulEl(D)$ then $1$-epimorphic\,?
\end{Question}
\eq

\section{Work of Cooke, Dwyer-Kan-Smith and Mitchell}
\label{SecCDKSM}

Let $G$ be a group, considered as a category. By a space we mean a topological space homotopy equivalent to a CW-complex. Let $(\text{Spaces})$ be the category of spaces and continuous maps.
Let $(\text{Hot})$ be the category of spaces and homotopy classes of continuous maps.

{\sc Cooke} developed in \bfcit{Co78} an obstruction theory for the induced functor
\[
\bo G,\, (\text{Spaces})\bc \;\lra\; \bo G,\, (\text{Hot})\bc
\]
to be dense. The obstructions are classes in the cohomology groups of $G$ with certain coefficients in dimensions $\ge 3$; cf.\ \bfcite{Co78}{Th.\ 1.1}.%

{\sc Dwyer,} {\sc Kan} and {\sc Smith} generalised this obstruction theory in \bfcit{DKS89} from a group $G$ to an arbitrary category $D$ (and even to topological categories).
The obstruction to the density of the according functor then are classes in the Hochschild-Mitchell groups of $D$ with certain coefficients in dimensions $\ge 3$, and a problem ``involving
fundamental groupoids''; cf.\ \bfcite{DKS89}{3.5,\,3.6}.

{\sc Mitchell} has given in \bfcit{Mi72} the following criterion for the Hochschild-Mitchell cohomological dimension of a poset to be less than or equal to $2$.

Given $n\ge 2$, the {\it suspended $n$-crown} $\text{SC}_n$ is the poset defined as follows. As a set, let $\text{SC}_n := \{ u_i,\, v_i\; :\; i\in\Z/n\} \disj\{ s,\, t\}$ consist of $2n + 2$ elements.
The partial ordering on $\text{SC}_n$ is generated by 
\[
v_i < u_i\; , \;\; v_i < u_{i-1} \; ,\;\; u_i < t \;\; \text{and} \;\; s < v_i \;\;\; \text{for all $i\in\Z/n$} \; .
\]

Suppose given a finite poset $D$. According to \bfcite{Mi72}{Th.\ 35.7}, its Hochschild-Mitchell cohomology vanishes in dimensions $\ge 3$ for all $D$-bimodules as coefficients if and only if
neither (i) nor (ii) holds.

\begin{itemize}
\item[(i)] The poset $D$ contains an isomorphic copy of $\text{SC}_n$ as a full subposet for some $n\ge 3$.
\item[(ii)] The poset $D$ contains an isomorphic copy of $\text{SC}_2$ as a full subposet, and there is no $d\in D$ such that $v_0\le d$, $v_1\le d$, $d\le u_0$ and $d\le u_1$.
\end{itemize}

\begin{Question}
\label{QuMit0}
Suppose a finite poset satisfies condition {\rm (i)} or {\rm (ii)}. Is it not ind-flat?
\end{Question}

I do not know a counterexample. An affirmative answer would provide a hint at a possible existence of obstruction classes in certain Hochschild-Mitchell cohomology groups in dimension $\ge 3$ against
the density of 
$$
\bo D, \El\bc \;\lra\; \bo D,\ulEl\bc
\leqno (\ast)
$$
for an arbitrary Frobenius category $\El$. For then ind-flat finite posets would have vanishing Hochschild-Mitchell cohomology in dimension $\ge 3$, so that, provided such an obstruction theory exists,
we would see the ``real reason'' why an ind-flat finite poset $D$ yields a dense functor $(\ast)$; cf.\ Theorem \ref{Th3.1.1}.

\bq
The following two simple examples should point out problems one possibly encounters when trying to approach Question \ref{QuMit0}.

\begin{Example}
\label{ExMit1}\rm
Let 
\[
D \; :=\; \big\{ \leer,\, \{ 1\},\,\{ 2\},\, \{ 3\},\,  \{ 1,2\},\, \{ 2,3\},\, \{ 1,3\},\, \{1,2,3\},\, \{ 1,3,4\},\, \{1,2,3,4\} \big\} \; ,
\]
ordered by inclusion. 
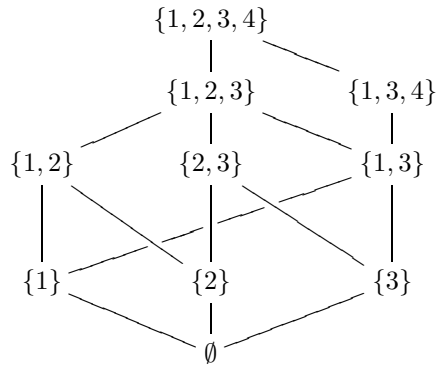
\[
\xymatrix@R=4mm{
                                & \{ 1,2,3,4\}                          &                                                \\
                                & \{ 1,2,3\}\ar@{-}[u]                  & \{1,3,4\}\ar@{-}[ul]                           \\
\{1,2\} \ar@{-}[ur]             & \{ 2,3\}\ar@{-}[u]                    & \{1,3\}\ar@{-}[ul]\ar@{-}[u]                   \\  
                                &                                       &                                                \\  
\{1\}\ar@{-}[uu]\ar@{-}[uurr]|(0.324){\phantom{\rule{2mm}{2mm}}}|(0.483){\phantom{\rule{2mm}{2mm}}}|(0.658){\phantom{\rule{2mm}{2mm}}} 
                                & \{2\}\ar@{-}[uu]\ar@{-}[uul]          & \{ 3\}\ar@{-}[uu]\ar@{-}[uul]                  \\  
                                & \leer\ar@{-}[ul]\ar@{-}[u]\ar@{-}[ur] &                                                \\  
}
\]
It contains the suspended $3$-crown \mb{$D\ohne\big\{\{1,2,3\},\, \{1,3\}\big\}$} as a full subposet. Moreover, $\{ 1,2,3,4\}$ is
a minimal element in \mb{$\VV(\{ 1,2\})\cap\VV(\{ 2,3\})\cap\VV(\{ 1,3,4\})$}. However, $\ind-crown\!\big(\Lambda^0(\{1,2,3,4\})\big)$ is componentwise $1$-connected. Only 
$\ind-crown\!\big(\Lambda^0(\{1,2,3\})\big)$ is not. 

Thus, if we are given a finite poset that contains a suspended $3$-crown with maximal element $t$, and even if, moreover, $t$ is chosen to be minimal with respect to lying over the rest of the 
suspended $3$-crown, we can still not conclude that $\ind-crown\!\big(\Lambda^0(t)\big)$ is not componentwise $1$-connected. Instead, we will have to search elsewhere for a suitable element $t'$ such 
that $\ind-crown\!\big(\Lambda^0(t')\big)$ is not componentwise $1$-connected in order to prove failure of ind-flatness.
\end{Example}

\begin{Example}
\label{ExMit2}\rm
Let 
\[
D \; :=\; \big\{ \leer,\, \{ 1\},\,\{ 2\},\, \{ 3\},\,  \{ 1,2\},\, \{ 2,3,4\},\, \{ 1,3,4\},\, \{ 1,2,3\},\,\{1,2,3,4\} \big\} \; ,
\]
ordered by inclusion. 
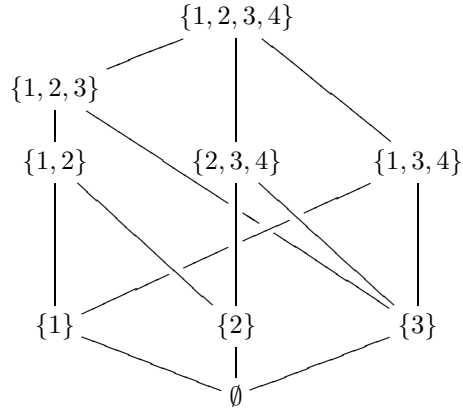
\[
\xymatrix@R=4mm{
                                & \{ 1,2,3,4\}                          &                                                \\
\{ 1,2,3\}\ar@{-}[ur]           &                                       &                                                \\
\{1,2\} \ar@{-}[u]              & \{ 2,3,4\}\ar@{-}[uu]                 & \{1,3,4\}\ar@{-}[uul]                          \\  
                                &                                       &                                                \\  
                                &                                       &                                                \\  
\{1\}\ar@{-}[uuu]\ar@{-}[uuurr]|(0.336){\phantom{\rule{2mm}{2mm}}}|(0.5){\phantom{\rule{2mm}{2mm}}}|(0.668){\phantom{\rule{2mm}{2mm}}} 
                                & \{2\}\ar@{-}[uuu]\ar@{-}[uuul]        & \{ 3\}\ar@{-}[uuu]\ar@{-}[uuul]\ar@{-}[uuuull]|(0.415){\phantom{\rule{2mm}{2mm}}}|(0.5){\phantom{\rule{2mm}{2mm}}} \\  
                                & \leer\ar@{-}[ul]\ar@{-}[u]\ar@{-}[ur] &                                                \\  
}
\]
It contains the suspended $3$-crown \mb{$D' := D\ohne\big\{\{1,2,3\}\big\}$} as a full subposet.
Now $\ind-crown\!\big(\Lambda_{D'}^0(\{ 1,2,3,4\})\big)$ is homotopy equivalent to a circle, whereas $\ind-crown\!\big(\Lambda_D^0(\{ 1,2,3,4\})\big)$ is 
homotopy equivalent to a wedge of two circles. So $D$ is not ind-flat. The reason for this, however, is an ind-crown of a somewhat surprising shape.

Thus if we want to attach some kind of homotopical invariant to a poset, or to a pair consisting of a poset and an element of it, and if we want to prove that this invariant is preserved under certain 
full embeddings of posets, we are confronted with this ``jump phenomenon''.

\end{Example}

\eq


\parskip0.0ex
\begin{footnotesize}

\parskip1.2ex

\begin{flushright}
Matthias K\"unzer\\
Lehrstuhl D f\"ur Mathematik\\
RWTH Aachen\\
Templergraben 64\\
D-52062 Aachen \\
kuenzer@math.rwth-aachen.de \\
www.math.rwth-aachen.de/$\sim$kuenzer\\
\end{flushright}
\end{footnotesize}

\end{document}